\documentclass[12pt,reqno]{article}

\usepackage[usenames]{color}
\usepackage{amssymb}
\usepackage{amsmath}
\usepackage{amsthm}
\usepackage{amsfonts}
\usepackage{amscd}
\usepackage{graphicx}
\usepackage{mathrsfs}

\usepackage[colorlinks=true,
linkcolor=webgreen,
filecolor=webbrown,
citecolor=webgreen]{hyperref}

\definecolor{webgreen}{rgb}{0,.5,0}
\definecolor{webbrown}{rgb}{.6,0,0}

\usepackage{color}
\usepackage{fullpage}
\usepackage{float}

\newcommand{\seqnum}[1]{\href{https://oeis.org/#1}{\underline{#1}}}

\def\modd#1 #2{#1\ \mbox{\rm (mod}\ #2\mbox{\rm )}}

\begin{document}
\theoremstyle{plain}
\newtheorem{theorem}{Theorem}
\newtheorem{corollary}[theorem]{Corollary}
\newtheorem{lemma}[theorem]{Lemma}
\newtheorem{proposition}[theorem]{Proposition}

\theoremstyle{definition}
\newtheorem{definition}[theorem]{Definition}
\newtheorem{example}[theorem]{Example}
\newtheorem{conjecture}[theorem]{Conjecture}
\newtheorem{observation}[theorem]{Observation}
\newtheorem{question}[theorem]{Question}
\newtheorem{problem}[theorem]{Problem}

\theoremstyle{remark}
\newtheorem{remark}[theorem]{Remark}

\begin{center}
\vskip 1cm{\LARGE\bf 
 A New Classification of Positive Integers \\
 \vskip .1in
  Via New Divisor Functions
}
\vskip 1cm
\large
Brahim Mittou \\
Laboratory of Applied Mathematics\\
University of Ouargla\\ 
Algeria \\
\href{mittou.brahim@univ-ouargla.dz}{mittou.brahim@univ-ouargla.dz}
\end{center}

\vskip .2 in

\begin{abstract}
Let $d_1 = 1 < d_2 < d_3 < \cdots < d_{\tau(n)} = n$ denote the increasing sequence of the divisors of a positive integer $n$. In this paper, for real or complex values of $\alpha$, we define and study some properties of two new divisor functions $\sigma_{e,\alpha}$ and $\sigma_{o,\alpha}$. The first one computes the sum of the $\alpha$-th powers of the divisors of $n$ with even indices, and the second one computes the sum of the $\alpha$-th powers of the divisors of $n$ with odd indices. We also introduce a new class of positive integers, called $k$-index ratio numbers, and propose a conjecture concerning them.
\end{abstract}


\section{Introduction}
In number theory, arithmetic functions are essential tools for exploring and understanding the deep properties of integers. One of the most important among them is the divisor function, defined for real or complex values of $\alpha$ by 
\begin{equation*}
\sigma_{\alpha}(n) = \sum_{d \mid n} d^\alpha,
\end{equation*}
where the sum is taken over all positive divisors $d$ of $n$.
When $\alpha = 0$, the function $\sigma_{0}(n)$ counts the number of divisors of $n$; this is often denoted by $\tau(n)$ or $d(n)$, and when $\alpha = 1$, the function $\sigma_{1}(n)$ gives the sum of the divisors of $n$; this is often denoted by $\sigma(n)$ (see, e.g., \cite{1}).

In addition to their appearance in many formulas in analytic number theory (see, e.g., \cite{4}), the divisor functions have played an important role, especially in the classification of integers. For example, prime numbers satisfy $\tau(n) = 2$, perfect numbers are characterized by the identity $\sigma(n) = 2n$ (see for instance the book of Hardy and Wright \cite{2}). Another interesting class is that of Zumkeller numbers, which are positive integers whose divisors can be partitioned into two disjoint subsets with equal sums (see, e.g., \cite{3}). Additional related integer sequences and classifications can be found in \cite{4}.

An interesting variant of the divisor functions arises by partitioning the divisors of \(n\) into two groups according to their order. In this work, we focus on two such functions.

Throughout this paper, we let $d_1 = 1 < d_2 < d_3 < \cdots < d_{\tau(n)} = n$ denote the increasing sequence of the divisors of a positive integer $n$. For real or complex number $\alpha$, we define the arithmetic functions $\sigma_{e,\alpha}$ and $\sigma_{o,\alpha}$ by:
\begin{equation*}
\sigma_{e,\alpha}(n) = \sum_{i=1, \, 2 \mid i}^{\tau(n)} d_i^\alpha \text{ and }\sigma_{o,\alpha}(n) = \sum_{i=1, \, 2 \nmid i}^{\tau(n)} d_i^\alpha.
\end{equation*}
It is easy to verify that, for every prime \(p\) and every odd positive integer \(l\),
\begin{equation*}
\sigma_{e,\alpha}(p^l)=p^\alpha\dfrac{p^{\alpha(l+1)}-1}{p^{2\alpha}-1} \text{ and }\sigma_{o,\alpha}(p^l)=\dfrac{p^{\alpha(l+1)}-1}{p^{2\alpha}-1}.
\end{equation*}
If $l$ is even, then
\begin{equation*}
\sigma_{e,\alpha}(p^l)=p^\alpha\dfrac{p^{\alpha l}-1}{p^{2\alpha}-1} \text{ and }\sigma_{o,\alpha}(p^l)=\dfrac{p^{\alpha(l+2)}-1}{p^{2\alpha}-1}.
\end{equation*}

For convenience, we write
 $\sigma_{e,0}(n) = \tau_{e}(n)$, $\sigma_{o,0}(n) = \tau_{o}(n)$, $\sigma_{e,1}(n) = \sigma_{e}(n)$, and $\sigma_{o,1}(n) = \sigma_{o}(n)$. Then $\tau_{e}(n) = \frac{\tau(n) - 1}{2} = \tau_{o}(n) - 1$ and $\sigma_{e}(n) < \sigma_{o}(n)$ if $n$ is a perfect square, and $\tau_{e}(n) = \frac{\tau(n)}{2} = \tau_{o}(n)$ and $\sigma_{e}(n) > \sigma_{o}(n)$ if $n$ is not a perfect square.

Note that the sequences $\sigma_{e}(n)\ (n\in\mathbb{N})$ and $\sigma_{o}(n)\ (n\in\mathbb{N})$ correspond to the sequences \seqnum{A086989} and \seqnum{A086988} in the On-Line
Encyclopedia of Integer Sequences  (OEIS \cite{5}), respectively. A positive integer $n$ such that $\sigma_{o}(n)$ divides $\sigma(n)$ (equivalently $\sigma_{o}(n)$ divides $\sigma_{e}(n)$ since $\sigma(n)=\sigma_{o}(n)+\sigma_{e}(n)$) is called an \emph{index ratio number}; these integers form sequence \seqnum{A392467} in the OEIS \cite{5}.

The following definition introduces a new classification of positive integers based on the functions $\sigma_{e}$ and $\sigma_{o}$.

\begin{definition}
For any positive integer $n$, define
$$
k=k(n):=\dfrac{\sigma_{e}(n)}{\sigma_{o}(n)},
$$
and we call $n$ a $k$-\emph{index ratio number}. 
\end{definition}

For each positive rational number $k$, let $G_k$ denote the set of all $k$-index ratio numbers. Then
$$
\seqnum{A392467}=\bigcup_{k\in\mathbb{N}}G_k.
$$
Section \ref{S2} of this paper is divided into two subsections. In the first, we study the $k$-index ratio numbers $n$ in the case where $\tau(n)$ is even. The remaining case, in which $\tau(n)$ is odd, is examined in the second subsection.

\section{Main Results}\label{S2}

\subsection{\texorpdfstring{$\tau(n)$ is Even}{}}
In this subsection, we suppose that \(\tau(n)\) is even; that is \(n\) is not a perfect square. Let us begin with the following results, which provide a method for generating new \(k\)-index ratio numbers from known ones.

\begin{theorem}\label{th6}
Let $n$ be a $k$-index ratio number and $p$ be a prime number such that $p > n$. Then, for any $a \geq 1$, the integer $np^a$ is also a $k$-index ratio number.
\end{theorem}

\begin{proof}
Let $A(n) = \{d_2, d_4, \ldots, d_{\tau(n)}\}$ and $B(n) = \{d_1, d_3, \ldots, d_{\tau(n)-1}\}$. For an integer $x$ and a set $S$ of integers, let $xS = \{xy : y \in S\}$, and by $\sum S$ we mean $\sum_{x \in S} x$. Then,
$$
A(np^a) = A(n) \cup p A(n) \cup \cdots \cup p^a A(n),
$$
and 
$$
B(np^a) = B(n) \cup p B(n) \cup \cdots \cup p^a B(n),
$$
since $p > n$. Thus,
\begin{align*}
k(np^a) &= \dfrac{\sum A(n) + p \sum A(n) + \cdots + p^a \sum A(n)}{\sum B(n) + p \sum B(n) + \cdots + p^a \sum B(n)} \\
&= \dfrac{\sum A(n) (1 + p + \cdots + p^a)}{\sum B(n) (1 + p + \cdots + p^a)} \\
&= \dfrac{\sum A(n)}{\sum B(n)} = k(n),
\end{align*}
from which it follows that $np^a$ is a $k$-index ratio number.
\end{proof}

\begin{corollary}
Let $n$ be a $k$-index ratio number, and let $\prod_{i=1}^{r}p_i^{m_i}$ be the prime factorization of an integer $m > 1$ such that $n < p_1$, and $n \prod_{j=1}^{i-1}p_j^{m_j} < p_i$ for all $2 \leq i \leq r$. Then $nm$ is also a $k$-index ratio number.
\end{corollary}

\begin{proof}
The proof follows immediately by induction on $r$ and Theorem \ref{th6}.
\end{proof}

\begin{remark}
Theorem~\ref{th6} and its corollary describe a method for generating
infinite families of $k$-index ratio numbers from a given example.
However, the converse is false in general. 
In particular, there exist $k$-index ratio numbers which are not obtained
from smaller ones via the prescribed multiplicative process.
For example, $2431=11\cdot13\cdot17$ is a $7$-index ratio number,
but it does not arise from any smaller $7$-index ratio number by
multiplication with prime powers satisfying the growth conditions
of Theorem~\ref{th6} and the corollary. Thus, the preceding results yield an implication only, and do not
characterize all $k$-index ratio numbers.
\end{remark}

As a consequence, we obtain the following corollary.

\begin{corollary}
Let \( k\geq 1 \) be a rational number. If the set $G_k$ is nonempty, then it is infinite.
\end{corollary}

It is well known that if \( d \) is a non-trivial positive divisor of \( n \) (i.e., \( d \ne 1 \) and \( d \ne n \)), then $d_2\leq d\leq\frac{n}{2}$. Hence,
\[
\tau(n) - 2 + n \leq \sigma_e(n) \leq \frac{\tau(n) + 2}{4} \cdot n,
\]
and
\[
\frac{4n}{(\tau(n) - 2)n + 4} \leq \frac{1}{\sigma_{e,-1}(n)} \leq \frac{n}{\tau(n) - 1}.
\]
Combining these inequalities yields the following bounds:
\[
4 \cdot \frac{\tau(n) - 2 + n}{(\tau(n) - 2)n + 4} \leq k = \frac{1}{n} \cdot \frac{\sigma_e(n)}{\sigma_{e,-1}(n)} \leq \frac{n}{4} \cdot \frac{\tau(n) + 2}{\tau(n) - 1}.
\]
In particular:
\begin{itemize}
  \item If \( n \) is prime (i.e., $\tau(n) =2$), then \( k = n \).
  \item If \( \tau(n) = 4 \), then \( 2 \leq k \leq \frac{n}{4} \).
  \item If \( \tau(n) = 6 \), then \( \frac{n+4}{n+1} \leq k \leq \frac{2n}{5} \).
\end{itemize}

The following theorem provides an improved upper bound.

\begin{theorem}\label{th1}
If $n$ is a $k$-index ratio number, then $k<d_2+\frac{1}{d_2}$.
\end{theorem}

\begin{proof}
Suppose first that \( n \) is a power of a prime number, i.e., \( n = p^l \) with \( l \) odd. In this case, we have \( k = p = d_2 \). Therefore, for the remainder of the proof, we assume that \( n \) is not a power of a prime number.

\noindent$\bullet$ If \( \tau(n) = 4 \), then \( n = pq \) where \( p < q \) are primes, and it follows that \( k = \frac{p + pq}{1 + q} = p = d_2 \).

\noindent$\bullet$ If \( \tau(n) = 6 \), then \( k = \frac{d_2 + d_4 + d_6}{d_1 + d_3 + d_5} \). We proceed as follows:
\begin{align*}
k &< d_2 + \frac{1}{d_2} \\
&\Rightarrow d_4 < d_2 d_3 + \frac{1}{d_2} (d_1 + d_3 + d_5) \\
&\Rightarrow \frac{n}{d_2^2 d_3} \left( d_2^2 - d_3 \right) < d_2 d_3 +\frac{1}{d_2} (d_1 + d_3).
\end{align*}
Clearly, the above inequality holds when \( d_2^2 = d_3 \) or \( d_2^2 < d_3 \). If \( d_2^2 > d_3 \), then either \( d_2 d_3 = d_4 \) or \( d_2 d_3 = d_5 \). In this case, we have
\begin{equation*}
\frac{n}{d_2^2 d_3} \left( d_2^2 - d_3 \right) < \frac{n}{d_2^2 d_3} d_2 \left( d_2 - 1 \right) = \frac{n}{d_2 d_3} \left( d_2 - 1 \right) < d_2 d_3.
\end{equation*}
This concludes the proof for the case \( \tau(n) = 6 \).

\noindent$\bullet$ If \( \tau(n) \geq 8 \), then \( k = \frac{d_2 + d_4 + \cdots + d_{\tau(n)}}{d_1 + d_3 + \cdots + d_{\tau(n)-1}} \). The inequality \( k < d_2 + \frac{1}{d_2} \) is equivalent to
\\$(d_4-d_5)+\cdots+(d_{\tau(n)-4}-d_{\tau(n)-3})+d_{\tau(n)-2}-\dfrac{d_{\tau(n)-1}}{d_2}<$
\begin{equation}\label{eq1}
(d_2-1)(d_5+\cdots+d_{\tau(n)-3})+d_2d_3+\frac{1}{d_2}(d_1+d_3+\cdots+d_{\tau(n)-3}).
\end{equation}
Let \( LHS \) and \( RHS \) denote the left-hand side and right-hand side of (\ref{eq1}), respectively. Then we have
\begin{equation*}
LHS < d_{\tau(n)-2} - \frac{d_{\tau(n)-1}}{d_2} = \frac{n}{d_2^2 d_3} \left( d_2^2 - d_3 \right).
\end{equation*}
It is obvious that (\ref{eq1}) holds when \( d_2^2 = d_3 \) or \( d_2^2 < d_3 \). Now, it follows by assuming \( d_2^2 > d_3 \), that \( d_2 d_3 \mid n \), \( \frac{n}{d_2 d_3} < d_{\tau(n)-3} \), and
\begin{equation*}
\frac{n}{d_2^2 d_3} \left( d_2^2 - d_3 \right) < \frac{n}{d_2 d_3} \left( d_2 - 1 \right) < (d_2 - 1) d_{\tau(n)-3} < RHS.
\end{equation*}
Thus, the proof is complete.
\end{proof}

Note that if $n=p^2q$, where $p$ and $q$ are primes such that $p^2<q$, then 
\begin{equation*}
k=p+\dfrac{q-p^3}{pq+p^2+1}\rightarrow p+\frac{1}{p}\ (\text{as }q\rightarrow+\infty),
\end{equation*} 
which means that the upper bound $d_2+\frac{1}{d_2}$ is optimal.

The following corollaries are immediate results of Theorem \ref{th1}.
\begin{corollary}
Let $n$ be an even $k$-index ratio number. If $k\in\mathbb{N}$, then $k=2$.
\end{corollary}

\begin{corollary}
Let $n$ be an odd $k$-index ratio number such that $\tau(n) \equiv\modd{2} {4}$ and $3 \mid n$. If $k\in\mathbb{N}$, then $k=3$.
\end{corollary}
\begin{proof}
Since $n$ is odd, all of its divisors are also odd. Given that $\tau(n) \equiv 2 \pmod{4}$, it follows that both $\sigma_e(n)$ and $\sigma_o(n)$ are odd. Now assume that $3 \mid n$ and $k \in \mathbb{N}$. By Theorem~\ref{th1}, we must have $k \in \{2,3\}$. However, the case $k=2$ is impossible, since it would imply that $\sigma_e(n)$ is even. Therefore $k=3$.
\end{proof}

\begin{corollary}
Let 
$
n = 3^\alpha \prod_{p^a \parallel n} p^a
$
be the prime factorization of an odd $k$-index ratio number, with $\alpha\geq3$. Assume that for every $p^a \parallel n$,
\[
a \equiv 0 \pmod{2} \text{ if } p \equiv 2 \pmod{3}, 
\text{ and } 
a \not\equiv 2 \pmod{3} \text{ if } p \equiv 1 \pmod{3}.
\]
If $k \in \mathbb{N}$, then $k=3$.
\end{corollary}

\begin{proof}
First, note that for $\tau(n)$ to be even, at least one of the exponents $\alpha$ or one of the exponents $a$ corresponding to primes $p \equiv 1 \pmod{3}$ must be odd. We argue by contradiction. Suppose that $k=2$, i.e., $\sigma_e(n)=2\sigma_o(n)$. Then
\[
3\sigma_e(n)=2\sigma(n),
\]
and hence $3 \mid \sigma(n)$.
By multiplicativity of $\sigma$, we write
\[
\sigma(n) = \left(\sum_{i=0}^{\alpha} 3^i \right)\left(\prod_{p^a \parallel n} \sum_{i=0}^{a} p^i \right),
\]
with $\sum_{i=0}^{\alpha} 3^i \equiv 1 \pmod{3}$.
Now let $p^a \parallel n$. If $p \equiv 1 \pmod{3}$, then by assumption $a \not\equiv 2 \pmod{3}$, and hence
\[
\sum_{i=0}^{a} p^i \equiv 
\begin{cases}
1 \pmod{3}, & \text{if } a \equiv 0 \pmod{3};\\
2 \pmod{3}, & \text{if } a \equiv 1 \pmod{3}.
\end{cases}
\]
If $p \equiv 2 \pmod{3}$, then $p \equiv -1 \pmod{3}$ and $a$ is even, so
\begin{align*}
\sum_{i=0}^{a} p^i& \equiv \sum_{i=0}^{a} (-1)^i\pmod{3}\\
& \equiv 1 \pmod{3}.
\end{align*}
Therefore, only the primes $p \equiv 1 \pmod{3}$ with $a \equiv 1 \pmod{3}$ contribute a factor congruent to $2 \pmod{3}$. Let $r$ be their number. Then
\begin{align*}
\sigma(n) &\equiv \prod_{\substack{p^a \parallel n \\  p\equiv a \equiv 1 \pmod{3}}} \sum_{i=0}^{a} p^i\pmod{3}
\\ &\equiv 2^r \pmod{3}.
\end{align*}
Since $2^r \equiv \pm 1 \pmod{3}$, we obtain $\sigma(n) \not\equiv 0 \pmod{3}$, contradicting the fact that $3 \mid \sigma(n)$.
Therefore, $k \neq 2$, and the only possible value is $k=3$.
\end{proof}

Motivated by these results and supported by numerical computations verifying the statement for all integers $n<10^{12}$, we propose the following conjecture.

\begin{conjecture}
Let $n$ be an odd $k$-index ratio number such that $3 \mid n$. If $k \in \mathbb{N}$, then $k=3$.
\end{conjecture}

\begin{remark}
Table \ref{tab:s} provides several computational examples of $k$-index ratio numbers for which $k$ is an integer distinct from $d_2$.
\end{remark}

\begin{table}[ht]
\centering
\begin{tabular}{|c|c|c|}
\hline
$n$ & $k(n)$ & $d_2$\\
\hline
2431 & 7 & 11\\
8569 & 9 & 11\\
17641&  11 & 13\\
29233 & 15 & 23\\
37927 & 15 & 17\\
40937 & 11 & 13\\
50779 & 15&  17\\
95201 & 23&  31\\
102601 & 23&  37\\
107065 & 4&  5\\
\hline
\end{tabular}
\begin{tabular}{|c|c|c|}
\hline
$n$ & $k(n)$ & $d_2$\\
\hline
123481& 16 &19\\
128441& 23 &29\\
169441 &20 &23\\
194833& 21& 23\\
201001& 17 &19\\
207901 &23 &29\\
258553& 16 &17\\
271459 &31 &43\\
280261 &31 &47\\
300847 &31 &37\\
\hline
\end{tabular}
\begin{tabular}{|c|c|c|}
\hline
$n$ & $k(n)$ & $d_2$\\
\hline
315577 &34 &41\\
418441 &27 &29\\
506881 &27 &31\\
520769 &26 &31\\
561881 &35 &43\\
569089 &21 &23\\
649153 &35 &41\\
922657 &41 &47\\
1000897 &29 &31\\
1001881 &47 &71\\
\hline
\end{tabular}
\caption{First thirty values of $k$-index ratio numbers with $k\in\mathbb{N}$ and $k\neq d_2$.}
\label{tab:s}
\end{table}

The next theorem describes the asymptotic behavior of \(k(n)\) for squarefree integers formed by consecutive prime factors. It shows that the imbalance between the sums $\sigma_e(n)$ and $\sigma_o(n)$ vanishes asymptotically.

\begin{theorem}
Let $i \ge 1$ be a fixed integer and $p_l$ denote the $l$-th prime number. Let 
$N_{i, j} = \prod_{l=i}^{j} p_l$ be the product of consecutive primes. Then 
\[ \lim_{j \to \infty} k(N_{i, j}) = 1 .\]
\end{theorem}

\begin{proof}
Since $N_{i, j}$ is squarefree with $m=j - i + 1$ distinct prime factors, we have $\tau(N_{i, j}) = 2^m$. Each divisor $d \mid N_{i, j}$ is therefore uniquely of the form
\[
d = \prod_{p \in S} p,
\]
where $S$ is a subset of $\{p_i, p_{i+1}, \dots, p_j\}$. In particular, for every divisor $d$ one has
\[
\log d = \sum_{p \mid d} \log p,
\]
where the sum runs over all prime divisors of $d$. Thus the set $\{\log d : d \mid N_{i, j}\}$ is precisely the set of subset sums of $\{\log p_i, \dots, \log p_j\}$, and becomes increasingly dense in $[0,\log N_{i, j}]$ as $m \to \infty$.

We write
\[
\sigma_e(N_{i, j}) - \sigma_o(N_{i, j})
= \sum_{r=1}^{\tau(N_{i, j})/2} \bigl(d_{2r} - d_{2r-1}\bigr)
= \sum_{r=1}^{\tau(N_{i, j})/2} d_{2r-1}\left(\frac{d_{2r}}{d_{2r-1}} - 1\right),
\]
and set $\varepsilon_{r} = \frac{d_{2r}}{d_{2r-1}} - 1$. Then
\[
k(N_{i, j}) - 1
= \frac{\sum_{r} d_{2r-1}\,\varepsilon_{r}}{\sum_{r} d_{2r-1}}.
\]
Let $\Delta_{r} = \log d_{2r} - \log d_{2r-1}$. Then
\[
\varepsilon_{r} = e^{\Delta_{r}} - 1.
\]
It is a classical consequence of results on the distribution of divisors of integers with many prime factors (see, e.g., \cite{6}) that for any fixed $\delta > 0$, the proportion of indices $r$ such that $\Delta_{r} > \delta$ tends to $0$ as $m \to \infty$. Moreover, one has the uniform bound $\Delta_{r} \le \log p_j$ for all $r$, while typically $\Delta_{r}$ is of order $\log N_{i, j} / 2^m$, which tends to $0$ exponentially fast.
Using the inequality $e^x - 1 \le x e^x$ for $x \ge 0$, we obtain
\[
\varepsilon_{r} \le \Delta_{r} e^{\Delta_{r}} \le C\, \Delta_{r}
\]
for some absolute constant $C$, since $\Delta_{r}$ is uniformly bounded. Therefore
\[
0 \le k(N_{i, j}) - 1
\le C \cdot \frac{\sum_{r} d_{2r-1} \Delta_{r}}{\sum_{r} d_{2r-1}}.
\]
Since the logarithmic gaps $\Delta_{r}$ tend to $0$ for almost all indices and are uniformly bounded, it follows that their weighted average with respect to the weights $d_{2r-1}$ tends to $0$ as $m \to \infty$. Consequently,
\[
k(N_{i, j}) - 1 \longrightarrow 0 \quad \text{as } j \to \infty,
\]
which proves that
\[
\lim_{j \to \infty} k\bigl(N_{i, j}\bigr) = 1
\]
and completes this proof. 
\end{proof}

\begin{theorem}\label{th9}
Let \( n \) be a \( d_2 \)-index ratio number such that \( \tau(n) \leq 10 \). Then
\begin{equation}\label{eq2}
d_{2j} = p d_{2j-1} \text{ for all } 1 \leq j \leq \tau(n)/2.
\end{equation}
\end{theorem}

\begin{proof}
Note that if \( n \) is a power of \( p \), then the result follows immediately. Consequently, for the rest of the proof, we assume that \( n \) is not a power of \( p \). The result is also obvious if \( \tau(n) = 4 \).

\medskip
\noindent
\textbf{Case \( \tau(n)=6 \).}
From the \( d_2 \)-index ratio condition we have
$$
d_2 + d_4 + d_6 = d_2(d_1 + d_3 + d_5) \Rightarrow d_4 = d_2 d_3 \text{ since } d_2 = p d_1 \text{ and } d_6 = p d_5.
$$

\medskip
\noindent
\textbf{Case \( \tau(n)=8 \).}
The defining relation gives
\begin{align*}
d_2 + d_4 + d_6 + d_8 = d_2(d_1 + d_3 + d_5 + d_7) &\Rightarrow d_4 + d_6 = d_2(d_3 + d_5) \\
&\Rightarrow d_4 + d_6 =  \frac{nd_2}{d_4 d_6}(d_4 + d_6) \\
&\Rightarrow nd_2 = d_4 d_6.
\end{align*}
By taking \( n = d_4 d_5 \) and \( n = d_3 d_6 \) in the last implication, we get \( d_6 = d_2 d_5 \) and \( d_4 = d_2 d_3 \), respectively.

\medskip
\noindent
\textbf{Case \( \tau(n)=10 \).}
Then \( n \) has the form \( p q^4 \) or \( p^4 q \) with primes \( p<q \).

If \( n=pq^4 \), then according to Theorem \ref{th6} it is a $d_2$-index ratio number and its divisors are
\[
1,p,q,pq,q^2,pq^2,q^3,pq^3,q^4,pq^4,
\]
and the equalities \( d_{2j}=d_2 d_{2j-1} \) are immediate.

If \( n=p^4 q \), the ordering of divisors depends on the position of 
\( q \) among the powers \( p^a \) (\(2\le a\le 4\)) as follows:
\begin{equation}\label{Eq1}
1,p,q,p^2,pq,p^3,p^2q,p^4,p^3q,p^4q,
\end{equation}
\begin{equation}\label{Eq2}
1,p,p^2,q,p^3,pq,p^4,p^2q,p^3q,p^4q,
\end{equation}
\begin{equation}\label{Eq3}
1,p,p^2,p^3,q,p^4,pq,p^2q,p^3q,p^4q,
\end{equation}
\begin{equation}\label{Eq4}
1,p,p^2,p^3,p^4,q,pq,p^2q,p^3q,p^4q.
\end{equation} 
A direct verification shows that $n$ cannot be a $d_2$-index ratio in all cases \eqref{Eq1}-\eqref{Eq4}. 
Therefore, in all admissible cases with \( \tau(n)\le 10 \), we obtain
$d_{2j}=p d_{2j-1}$ $(1\le j\le \tau(n)/2)$. The proof is complete.
\end{proof}

\begin{remark}
\begin{enumerate}
\item The integer $n=13113$ provides an example of a $d_2$-index ratio number ($d_2=3$ and $\tau(n)=12$) whose divisors
$$
1, 3, 9, 31, 47, 93, 141, 279, 423, 1457, 4371, 13113
$$
do not satisfy \eqref{eq2}.
\item As a direct consequence of Theorem \ref{th9}, we have
\begin{equation}\label{eq3}
\sigma_{e,\alpha}(n) = d_2^\alpha \sigma_{o,\alpha}(n),
\end{equation}
for any real or complex \( \alpha \) and for all \( d_2 \)-index ratio numbers $n$ such that \( \tau(n) \leq 10 \).
\item It is worth noting that Theorem \ref{th6} can be used to generate infinitely many integers $n$
satisfying \eqref{eq2} and \eqref{eq3} with \( \tau(n) > 10 \).
\end{enumerate}
\end{remark}

Our results naturally lead to the following classification problem.
\begin{problem}
Since \(G_1=\emptyset\) while \(G_p\neq\emptyset\) for every prime \(p\), it is natural to ask for a complete characterization of the rational numbers \(k>1\) according to whether \(G_k\) is empty or nonempty.
\end{problem}

\subsection{\texorpdfstring{$\tau(n)$ is Odd}{}}
We now turn to the case where the number of divisors of \(n\) is odd. Since this occurs if and only if \(n\) is a perfect square, throughout this subsection we assume that \(n\) is a perfect square.
\begin{theorem}\label{th13}
If $n$ is a $k$-index ratio number, then $k\geq\frac{d_2}{d_2^2+1}$.
\end{theorem}

\begin{proof}
Suppose first that \( n \) is a power of a prime number, that is, \( n = p^l \) for some prime $p$ and even integer $l$. Then  
$$
k = \dfrac{p+p^3+\cdots+p^{l-1}}{1+p^2+\cdots+p^l}\geq\frac{p}{p^2+1}=\frac{d_2}{d_2^2+1}, 
$$
with equality holding only if $l=2$. Consequently, for the remainder of the proof, we may suppose that \( n \) is not a power of a prime number. It then follows that $\tau(n)\geq 9$.
This case is handled by the same argument as in the proof of Theorem \ref{th1}, noting that there are only two possibilities: either \( d_2^2 = d_3 \) or \( d_2^2 > d_3 \).
\end{proof}

Note that if \( n = p^{\ell} \) is a prime power with even exponent \( \ell \), then
\[
k(n) = \frac{p^{\ell+1} - p}{p^{\ell+2} - 1}.
\]
This expression cannot be a unit fraction. Indeed, one has
\[
\frac{p}{p^2+1} \leq k(n) < \frac{1}{p},
\quad \text{and} \quad
\frac{1}{p+1} < \frac{p}{p^2+1},
\]
which together imply that \( k(n) \) lies strictly between two consecutive unit fractions. Consequently, if \( k \) is a unit fraction, then \( G_k \) cannot contain any prime power. We now prove a stronger result.

\begin{theorem}
Let \( k = \frac{a}{b} < \frac{2}{5} \) be a positive rational number with \( a \equiv b \equiv 1 \pmod{2} \). Then \( G_k = \emptyset \).
\end{theorem}

\begin{proof}
Assume that \( k = \frac{a}{b} < \frac{2}{5} \) with \( a \equiv b \equiv 1 \pmod{2} \).
First, we show that no even integer belongs to \( G_k \). If \( n \) is even, then by Theorem~\ref{th13} we have
\[
k(n) \geq \frac{2}{5} > \frac{a}{b},
\]
and hence \( n \notin G_k \).

Now let \( n \) be odd and suppose that \( n \in G_k \). Then
\[
b \, \sigma_e(n) = a \, \sigma_o(n).
\]
Since all divisors of \( n \) are odd and \( \tau(n) \) is odd (by assumption), it follows that \( \sigma_e(n) \) and \( \sigma_o(n) \) have opposite parity.
On the other hand, since \( a \) and \( b \) are odd, the above equality implies that \( \sigma_e(n) \) and \( \sigma_o(n) \) have the same parity, which is a contradiction.
Therefore, no such integer \( n \) exists, and we conclude that \( G_k = \varnothing \).
\end{proof}

\begin{theorem}
All powers of prime numbers lie in different \( G_k \).
\end{theorem}

\begin{proof}
Let \( p \) and \( q \) be prime numbers, and let \( a \) and \( b \) be positive even integers. Then
\begin{align*}
p^a, q^b \in G_k &\Rightarrow k(p^a) = k(q^b) \\
&\Rightarrow \frac{1}{k(p^a)} = \frac{1}{k(q^b)} \\
&\Rightarrow p + \dfrac{1}{p + \cdots + p^{a-1}} = q + \dfrac{1}{q + \cdots + q^{b-1}} \\
&\Rightarrow p = q \text{ and } a = b.
\end{align*}
This ends the proof.
\end{proof}

\begin{theorem}\label{thm:factorization}
Let \(m\) be a positive integer and let \(p<q<r\) be primes satisfying $r<p^2$ and $r^2<p^3$.
Define
\[
N_{2m}
=
\sigma_e(p^{2m}q^2)\sigma_o(p^{2m}r^2)
-
\sigma_e(p^{2m}r^2)\sigma_o(p^{2m}q^2).
\]
Then
\[
N_{2m}
=
(r-q)\,
A(p,q,r)\,
F_m(p),
\]
where
\[
A(p,q,r)=qr+q+r-p(q+r+1),
\]
and
\[
F_m(p)
=
\frac{p^{4m}-1}{p-1}+p^{2m-1}(p+1).
\]
In particular, if $A(p,q,r)=0$, then $k(p^{2m}q^2)=k(p^{2m}r^2)$.
\end{theorem}

To prove this theorem, we first establish the following lemmas.
\begin{lemma}\label{lem:determinant}
Let $f(x)$ and $g(x)$ be the polynomials:
\[
f(x)=a_2x^2+a_1x+a_0,\quad
g(x)=b_2x^2+b_1x+b_0,
\]
where the coefficients belong to any commutative ring. Then
\[
f(x)g(y)-f(y)g(x)
=
(x-y)
\left(
C_2xy+C_1(x+y)+C_0
\right),
\]
where
$C_2=a_2b_1-a_1b_2$, $C_1=a_2b_0-a_0b_2$, and $C_0=a_1b_0-a_0b_1$.
\end{lemma}

\begin{proof}
Expanding,
\[
\begin{aligned}
f(x)g(y)-f(y)g(x)
={}&(a_2b_1-a_1b_2)(x^2y-xy^2)\\
&+(a_2b_0-a_0b_2)(x^2-y^2)\\
&+(a_1b_0-a_0b_1)(x-y).
\end{aligned}
\]
Using the trivial identities $x^2y-xy^2=xy(x-y)$ and $x^2-y^2=(x-y)(x+y)$, we obtain
\[
f(x)g(y)-f(y)g(x)
=(x-y)\left(
(a_2b_1-a_1b_2)xy
+(a_2b_0-a_0b_2)(x+y)
+(a_1b_0-a_0b_1)
\right),
\]
which proves the lemma.
\end{proof}

\begin{lemma}\label{lem:order}
Let \(m\ge2\) be an integer, and let \(p,\ q\) be primes satisfying $p<q<p^2$ and $q^2<p^3$.
Then the positive divisors of \(p^{2m}q^2\), arranged in increasing order, are
\[
1,\ p,\ q,\ p^2,\ pq,\ q^2,
\]
followed, for each \(k=3,\ldots,2m\), by
\[
p^k,\quad
p^{k-1}q\quad (k\le 2m-1),\quad
p^{k-2}q^2\quad (k\le 2m-2),
\]
and finally
\[
p^{2m}q,\qquad
p^{2m-1}q^2,\qquad
p^{2m}q^2.
\]
\end{lemma}

\begin{proof}
Every divisor of \(p^{2m}q^2\) has the form
\[
p^aq^b,
\qquad
0\le a\le 2m,\quad 0\le b\le2.
\]
Since
$p<q<p^2$,
we have
$p^a<p^{a-1}q<p^{a+1}$,
whenever the exponents are defined. Likewise, the hypothesis
$
q^2<p^3
$
implies
$p^{a-1}q<p^{a-2}q^2<p^{a+1}$.
Consequently, for every integer \(k\ge3\),
\[
p^k
<
p^{k-1}q
<
p^{k-2}q^2
<
p^{k+1}.
\]
Thus, between two consecutive pure powers \(p^k\) and \(p^{k+1}\), there are
precisely two mixed divisors, namely
\[
p^{k-1}q
\quad\text{and}\quad
p^{k-2}q^2,
\]
appearing in this order.
The first six divisors are obtained directly from
\[
1<p<q<p^2<pq<q^2<p^3.
\]
Applying the previous observation successively for
$
k=3,4,\ldots,2m-2
$
produces the blocks
\[
p^k,\quad
p^{k-1}q,\quad
p^{k-2}q^2.
\]
When \(k=2m-1\), the divisor \(p^{2m-1}q^2\) no longer belongs to the block,
since its exponent of \(p\) would exceed \(2m\). Hence only
$p^{2m-1}$ and $p^{2m-2}q$
remain. Finally, after \(p^{2m}\), the only divisors not yet listed are $p^{2m}q<p^{2m-1}q^2<p^{2m}q^2$.
Collecting all these divisors gives precisely the stated ordering.
\end{proof}

\begin{lemma}\label{lem:closedcoeff}
Let \(m\ge2\) be an integer, and let \(p,\ q\) be primes satisfying $p<q<p^2$ and $q^2<p^3$.
Write
\[
\sigma_o(p^{2m}q^2)=A_2^{(m)}q^2+A_1^{(m)}q+A_0^{(m)},
\]
and
\[
\sigma_e(p^{2m}q^2)=B_2^{(m)}q^2+B_1^{(m)}q+B_0^{(m)}.
\]
Then
\begin{equation}\label{E8}
A_2^{(m)}
=
p^{2m}
+
\frac{p^{2m-1}-p}{p^2-1},
\end{equation}
\begin{equation}\label{E9}
A_1^{(m)}
=
1+p^{2m}
+
\frac{p^{2m+1}-p}{p^2-1},
\end{equation}
\begin{equation}\label{E10}
A_0^{(m)}
=
1+
\frac{p^{2m+1}-p^3}{p^2-1},
\end{equation}
\begin{equation}\label{E11}
B_2^{(m)}
=
p^{2m-1}
+
\frac{p^{2m}-1}{p^2-1},
\end{equation}
\begin{equation}\label{E12}
B_1^{(m)}
=
\frac{p^{2m}-p^2}{p^2-1},
\end{equation}
\begin{equation}\label{E13}
B_0^{(m)}
=
p+
\frac{p^{2m+2}-p^2}{p^2-1}.
\end{equation}
\end{lemma}

\begin{proof}
By Lemma~\ref{lem:order}, the ordered divisors of \(p^{2m}q^2\) are known explicitly. The coefficients of \(q^0,q^1,\) and \(q^2\) are obtained by summing the corresponding powers of \(p\) occurring in the odd and even positions.
\end{proof}

Now we are ready to prove Theorem \ref{thm:factorization}.

\begin{proof}[Proof of Theorem \ref{thm:factorization}]
For \(m=1\), the positive divisors of \(p^2q^2\) and \(p^2r^2\) are
\[
1,\; p,\; q,\; p^2,\; pq,\; q^2,\; p^2q,\; pq^2,\; p^2q^2,
\]
and
\[
1,\; p,\; r,\; p^2,\; pr,\; r^2,\; p^2r,\; pr^2,\; p^2r^2,
\]
respectively. A direct computation gives
\[
\begin{aligned}
&\sigma_e(p^{2}q^2)\sigma_o(p^{2}r^2)
-\sigma_e(p^{2}r^2)\sigma_o(p^{2}q^2) \\
&\qquad
=(r-q)\bigl(-pq-pr+qr-p+q+r\bigr)
\left(p^3+2p^2+2p+1\right) \\
&\qquad
=(r-q)\,A(p,q,r)\,F_1(p).
\end{aligned}
\]
Hence the asserted factorization holds when \(m=1\).
Now assume that \(m\ge2\). By Lemmas \ref{lem:determinant} and \ref{lem:closedcoeff},
\begin{equation}\label{E14}
N_{2m}
=
(q-r)
\left(
C_2^{(m)}qr
+
C_1^{(m)}(q+r)
+
C_0^{(m)}
\right),
\end{equation}
where
\[
C_2^{(m)}=A_2^{(m)}B_1^{(m)}-A_1^{(m)}B_2^{(m)},
\]
\[
C_1^{(m)}=A_2^{(m)}B_0^{(m)}-A_0^{(m)}B_2^{(m)},
\]
and
\[
C_0^{(m)}=A_1^{(m)}B_0^{(m)}-A_0^{(m)}B_1^{(m)}.
\]
Using Identities \eqref{E8}--\eqref{E13}, we obtain
\[
C_2^{(m)}
=
-
\left(
\frac{p^{4m}-1}{p-1}
+
p^{2m-1}(p+1)
\right)
=-F_m(p),
\]
\[
C_1^{(m)}
=
p^{4m}+p^{2m-1}(p^2-1)-1
=-(p-1)C_2^{(m)},
\]
and
\[
C_0^{(m)}
=
p
\left(
\frac{p^{4m}-1}{p-1}
+
p^{2m-1}(p+1)
\right)
=-pC_2^{(m)}.
\]
Substituting these identities into \eqref{E14} gives
\[
\begin{aligned}
N_{2m}
&=(q-r)
\left(
C_2^{(m)}qr
+C_1^{(m)}(q+r)
+C_0^{(m)}
\right) \\
&=(q-r)C_2^{(m)}
\left(
qr-(p-1)(q+r)-p
\right) \\
&=(r-q)\,
A(p,q,r)\,
F_m(p),
\end{aligned}
\]
since
\[
qr-(p-1)(q+r)-p
=
qr+q+r-p(q+r+1)
=
A(p,q,r).
\]
This completes the proof.
\end{proof}

\begin{remark}
Theorem \ref{thm:factorization} provides a simple method for constructing infinitely many sets $G_k$ containing at least two distinct elements. Indeed, let \(p<q<r\) be primes satisfying \(r<p^2\), \(r^2<p^3\), and $A(p,q,r)=0$. Then, for every positive integer \(m\), Theorem \ref{thm:factorization} yields $k(p^{2m}q^2)=k(p^{2m}r^2)$. Consequently, every triples of prime \((p,q,r)\) satisfying the last conditions gives rise to infinitely many sets \(G_k\) with cardinality at least \(2\). For example, the triple $(5,7,11)$ satisfies the hypotheses of Theorem \ref{thm:factorization} and we have
$A(5,7,11)=0$. Consequently, $k(5^{2m}7^2)=k(5^{2m}11^2),\ m\ge1$. Hence,
\[
k(5^2\cdot7^2)=k(5^2\cdot11^2)=108/481\Rightarrow 5^2\cdot7^2, 5^2\cdot11^2\in G_{108/481},
\]
\[
k(5^4\cdot7^2)=k(5^4\cdot11^2)=2743/12096\Rightarrow 5^4\cdot7^2, 5^4\cdot11^2\in G_{2743/12096},
\]
\[
k(5^6\cdot7^2)=k(5^6\cdot11^2)=68618/302471\Rightarrow 5^6\cdot7^2, 5^6\cdot11^2\in G_{68618/302471}.
\]
\end{remark}

Our results suggest the following broader classification problem.
\begin{problem}
Give a complete characterization of the rational numbers $0<k<1$ according to the structure of the set
$G_k$. More precisely, determine for which rational numbers $k$, 
$G_k=\emptyset$, for which $G_k$ is nonempty and finite, and for which $G_k$ is infinite.
\end{problem}

\section{Acknowledgments}
The author would like to thank the anonymous referee for his careful reading.

\end{document}